\documentclass{amsart}
\usepackage{amsmath}
\usepackage{amssymb}
\usepackage{amscd}
\usepackage{latexsym}
\usepackage{graphicx} 

\theoremstyle{plain}
\newtheorem{thm}{Theorem}
\newtheorem{lem}{Lemma}
\newtheorem{cor}{Corollary}

\theoremstyle{definition}

\def \CPb {\overline{\mathbf{CP}}{}^{2}}
\def \CP {\mathbf{CP}{}^{2}} 
\def \CPC {\mathbf{CP}{}^{2}\# \,3\overline{\mathbf{CP}}{}^{2}}
\def \R {\mathbf{R}}
\def \Z {\mathbf{Z}}
\def \Sig{\Sigma}

\def \SS {S^2\times S^2}

\def \vp {\varphi}

\def \a {\alpha}
\def \b {\beta}
\def \g {\gamma}

\def \lam {\lambda}
\def\L{\Lambda}

\def \o {\omega}

\def \z {\zeta}
\def \ba {\bar{\alpha}}
\def \bd {\partial}
\def \x {\times}
\def \- {\!\smallsetminus\!}
\def \C {\subset}

\def \ba {\bar{a}}

\def \hg {\hat{g}}
\def \ha {\hat{a}}

\def \eu{{\text{e}}}
\def \sign{{\text{sign}}}

\def \ssw {\text{SW}}

\def \swu{\text{SW}'}

\def\fs{\mathfrak{s}}

\def\hk{\widehat{k}}

\begin{document}

\baselineskip.5cm
\title {Reverse engineering small 4-manifolds}
\author[Ronald Fintushel]{Ronald Fintushel}
\address{Department of Mathematics, Michigan State University \newline
\hspace*{.375in}East Lansing, Michigan 48824}
\email{\rm{ronfint@math.msu.edu}}
\thanks{The first author was partially supported NSF Grant DMS0305818
and the third author by NSF Grant DMS0505080. The second author was partially supported by CFI, NSERC, and OIT grants.\\ \indent R.F. and R.J.S. also wish to thank the institutions they visited while the ideas in this paper percolated: Banff International Research Station, Mathematical Sciences Research Institute,  and The Centre for Mathematical Sciences, University of Cambridge. B.D.P. wishes to thank Anar Akhmedov for helpful conversations.}
\author[B. Doug Park]{B. Doug Park}
\address{Department of Pure Mathematics, University of Waterloo \newline
\hspace*{.375in}Waterloo, Ontario N2L 3G1, Canada}
\email{\rm{bdpark@math.uwaterloo.ca}}
\author[Ronald J. Stern]{Ronald J. Stern}
\address{Department of Mathematics, University of California \newline
\hspace*{.375in}Irvine,  California 92697}
\email{\rm{rstern@uci.edu}}
\begin{abstract} We introduce a general procedure called `reverse engineering'  that can be used to construct infinite families of smooth $4$-manifolds in a given homeomorphism type. As one of the applications of this technique, we produce an infinite family of pairwise nondiffeomorphic $4$-manifolds homeomorphic to $\CPC$.
\end{abstract}
\maketitle

\section{Introduction\label{Intro}}

This paper introduces a technique which we call  {\it reverse engineering}\/ that can be used to construct infinite families of distinct smooth structures on many $4$-manifolds. As one example of the utility of this technique we will construct infinitely many distinct smooth structures on $\CPC$. Exotic smooth structures on these manifolds were first constructed in \cite{AP, BK}. 

Reverse engineering is a three step process for constructing infinite families of distinct smooth structures on a given simply connected $4$-manifold.  One starts  with a model manifold which has nontrivial Seiberg-Witten invariant and the same euler number and signature as the simply connected manifold $X$ that one is trying to construct, but with $b_1>0$. The second step is to  find $b_{1}$ essential tori that carry generators of $H_1$ and to surger each of these tori in order to kill $H_1$ and, in favorable circumstances, to kill $\pi_{1}$. The third step is to compute Seiberg-Witten invariants. After each of the first  $b_{1}-1$ surgeries one needs to preserve the fact that the Seiberg-Witten invariant is nonzero. The fact that the next to last manifold in the string of surgeries has nontrivial Seiberg-Witten invariant allows the use of the Morgan, Mrowka, Szab\'{o} formula \cite{MMS} to produce an infinite family as was done in \cite{Santeria}. 

In many instances this procedure can be successfully applied without any computation, or even mention, of Seiberg-Witten invariants. If the model manifold for $X$ is symplectic and $b_{1}-1$ of the tori are Lagrangian so that a Luttinger surgery will reduce $b_{1}$, then there are infinitely many distinct smooth manifolds with the same cohomology ring as $X$. If the resulting manifold is simply connected, then one can often show that there are infinitely many distinct smooth structures on $X$. Aside from finding interesting model manifolds, it seems that the most difficult aspect to the reverse engineering procedure is the  computation of fundamental groups.

We will prove the main theorem that shows that this procedure provides infinitely many distinct manifolds in \S2. We then provide two examples. In \S3 we apply the reverse engineering  procedure to a model for $\CPC$,  the $2$-fold symmetric product of a genus 3 surface. We will identify the Lagrangian tori, show that the Luttinger surgeries result in a simply connected manifold, and produce infinitely many distinct smooth structures on $\CPC$. In \S4 we apply the reverse engineering  procedure  to the product of two genus $2$ surfaces, a model for $\SS$. We will identify Lagrangian tori that kill $H_{1}$ and construct infinitely many distinct smooth manifolds with the cohomology ring of $\SS$. We have been unsuccessful in showing that these  manifolds are simply connected.

\section{Reverse Engineering}

One of the key questions in smooth $4$-manifold topology is whether a fixed homeomorphism type containing a smooth $4$-manifold must actually contain infinitely many diffeomorphism types. The idea of this section is to state and prove a general theorem pointing in this direction which may be useful to those who are constructing exotic $4$-manifolds. 

To state our theorem, we need to discuss some notation related to surgery on a  torus with trivial normal bundle.  Suppose that $T$ is such a torus with tubular neighborhood $N_T$. Let $\a$ and $\b$ be generators of $\pi_1(T^2)$ and  let $S^1_{\a}$ and $S^1_{\b}$ be loops in $T^3=\bd N_T$ homologous in $N_T$ to $\a$ and $\b$ respectively. Let $\mu_T$ denote a meridional circle to $T$ in $X$. By {\it{$p/q$-surgery}} on $T$ with respect to $\b$ we mean 
\begin{gather*} X_{T,\b}(p/q) =(X\- N_T) \cup_{\vp} (S^1\x S^1 \x D^2),\\
\vp: S^1\x S^1 \x \bd D^2 \to \bd(X\- N_T) 
\end{gather*}
where the gluing map 
satisfies $\vp_*([\bd D^2]) = q [S^1_{\b}] + p [\mu_T]$ in $H_1(\bd(X\- N_T);\Z )$. We denote the `core torus' $S^1\x S^1\x\{0\}\C X_{T,\b}(p/q)$ by $T_{p/q}$. 

We have framed $N_{T}$ using $S^1_{\a}$ and $S^1_{\b}$; so the pushoffs of $\a$ and $\b$ in this framing are $S^1_{\a}$ and $S^{1}_{\b}$. When the curve $S^{1}_{\b}$ is nullhomologous in $X\- N_T$, then $H_1(X_{T,\b}(1/q);\Z )= H_1(X;\Z )$. In addition, when $T$ itself is nullhomologous,  then $H_1(X_{T,\b}(p/q);\Z )= H_1(X;\Z ) \oplus \Z / p\Z$.
 
If $X$ is a symplectic manifold and $T$ is  any Lagrangian torus, then there is a canonical framing, called the Lagrangian framing, of $N_{T}$. This framing is uniquely determined by the property that pushoffs of $T$ in this framing remain Lagrangian. If one performs $1/n$ surgeries with respect to the pushoff in this framing of any curve on $T$, then the result is also a symplectic manifold. We refer the reader to  \cite{ADK} for a full discussion of this phenomenon, which is referred to there as {\it{Luttinger surgery}}. One must be careful to note that if the pushoff of a curve using the Lagrangian framing is not nullhomologous in $X\- N_T$, then $1/n$ surgery may in fact change $H_{1}$. 

Our theorem is:

\begin{thm}\label{SW} Let $X$ be a smooth closed oriented $4$-manifold which contains a nullhomologous torus $\L$, and let $\lam$ be a simple loop on $\L$ so that $S^1_{\lam}$ is nullhomologous in $X\- N_{\lam}$. If the Seiberg-Witten invariant of $X_{\L,\lam}(0)$ is nontrivial in the sense that for some basic class $k_0$, $\sum\limits_i\swu_{X_{\L,\lam}(0)}(k_0+2i[\L_0])\ne0$, then among the manifolds $\{ X_{\L,\lam}(1/n)\}$, infinitely many are pairwise nondiffeomorphic.
\end{thm}

The meaning of `$\swu$' is explained below. The following is a very simple but effective corollary to the proof.

\begin{cor} \label{CSW} Suppose that $X_0=X_{\L,\lam}(0)$ has, up to sign, just one Seiberg-Witten basic class. Then the manifolds $X_n=X_{\L,\lam}(1/n),\ n=1,2,3,\dots$ are pairwise nondiffeomorphic. \qed
\end{cor}

In case $X$ is simply connected and the  $\{ X_{\L,\lam}(1/n)\}$ are also simply connected then all of the manifolds $\{ X_{\L,\lam}(1/n)\}$ are homeomorphic.

As outlined in the introduction, one very useful application of the theorem is to start with a model manifold  with $b_1>0$  and nontrivial Seiberg-Witten invariant and with the same euler number and signature as a (say) $b_1=0$ manifold that we are trying to construct. Then, provided that we can find them, we surger essential tori which carry generators of $H_1$. If we can do this $b_1$ times, we kill $b_1$. We will see that if we can arrange the Seiberg-Witten invariant to be nonzero after each surgery, then  
we will satisfy the hypothesis that $X_{\L,\lam}(0)$ have nontrivial Seiberg-Witten invariant. For example, if we start with a symplectic $4$-manifold with $b^+>1$ and each time perform a Luttinger surgery on an embedded Lagrangian torus, this will be true. The fact that the next to last manifold in our string of surgeries has nontrivial Seiberg-Witten invariant will allow the application of the theorem and/or its corollary. We will discuss two examples in \S3 and \S4.

The proof of Theorem~\ref{SW} involves calculation of Seiberg-Witten invariants. We give a short discussion for the purpose of setting notation. The Seiberg-Witten
invariant  of a smooth closed oriented $4$-manifold
$X$ with a homology orientation, {\it{i.e.}} an orientation of $(H^0\oplus H^1\oplus H^+)(X;\R)$, and with $b^+_X>1$ is an integer-valued function $\ssw_X$ which is defined on the set of $spin ^c$ structures over $X$. Corresponding to each $spin ^c$ structure $\fs $ over $X$ is the bundle of positive spinors $W^+_{\fs}$ over $X$. Set $c(\fs)\in H_2(X;\Z)$ to be the 
Poincar\'e dual of $c_1(W^+_{\fs})$.
Each $c(\fs)$ is a characteristic element of $H_2(X;\Z)$ (i.e. its
Poincar\'e dual $\hat{c}(\fs)=c_1(W^+_{\fs})$ reduces to $w_2(X)$ mod~2).
We shall work with the modified
Seiberg-Witten invariant 
\[ \swu_X: \lbrace k\in H_2(X;\Z)\mid\hat{k}\equiv w_2(X)\pmod2\rbrace
\rightarrow \Z \]
defined by $\swu_X(k)=\sum\limits_{c(\fs)=k}\ssw_X(\fs)$. This is a diffeomorphism  invariant of $X$ together with its homology orientation. (As usual, the choice of homology orientation is suppressed from notation.)
If $H_1(X;\Z)$ has no $2$-torsion, then $\swu_X=\ssw_X$.

In case $b^+_X=1$, the invariant requires the choice of a class $H\in H_2(X;\R)$ with 
$H\cdot H>0$. We now need to be a bit more explicit.  Suppose we have a given orientation of $H^2_+(X;\R)$ and a given metric for $X$. The Seiberg-Witten invariant depends on the metric $g$ and a self-dual 2-form as follows. 
There is a unique $g$-self-dual harmonic 2-form $\o_g\in H^2_+(X;\R)$ with $\o_g^2=1$ and corresponding to the positive orientation. Fix a characteristic homology class $k\in H_2(X;\Z)$.  Given a pair $(A,\psi)$, where
$A$ is a connection in the complex line bundle whose first Chern class is the Poincar\'e dual $\hk=\frac{i}{2\pi}[F_A]$ of $k$ and $\psi$ a section of the bundle of self-dual spinors for the associated $spin^{\,c}$ structure, the perturbed Seiberg-Witten equations are:
\[
D_A\psi = 0, \ \ 
F_A^+  = q(\psi)+i\eta \notag\label{SWeqn},
\]
where $F_A^+$ is the self-dual part of the curvature $F_A$,
$D_A$ is the twisted Dirac operator, $\eta$ is a
self-dual 2-form on $X$, and
$q$ is a quadratic function. Write $\ssw_{X,g,\eta}(k)$ for the
corresponding invariant. As the pair
$(g,\eta)$ varies, $\ssw_{X,g,\eta}(k)$ can change only at those pairs
$(g,\eta)$ for which there are solutions with $\psi=0$. These 
solutions occur for pairs $(g,\eta)$ satisfying $(2\pi\hk+\eta)\cdot\o_g=0$.
This last equation defines a wall in $H^2(X;\R)$. 

The point $\o_g$ determines a component of the double cone consisting of elements of $H^2(X;\R)$ of positive square. We prefer to work with $H_2(X;\R)$. The dual component is determined by the Poincar\'e dual $H$ of $\omega_g$.  An element $H'\in H_2(X;\R)$ of positive square lies in the same component as $H$ if $H'\cdot H>0$. If
$(2\pi \hk+\eta)\cdot\o_g\ne 0$ for a generic $\eta$, $\,\ssw_{X,g,\eta}(k)$ is
well-defined, and its value depends only on the sign of $(2\pi \hk+\eta)\cdot\o_g$. Write $\ssw_{X,H}^+(k)$ for $\ssw_{X,g,\eta}(k)$ if 
$(2\pi \hk+\eta)\cdot\o_g>0$ and $\ssw_{X,H}^-(k)$ in the other case.

The invariant $\ssw_{X,H}(k)$ is defined by $\ssw_{X,H}(k) =\ssw_{X,H}^+(k)$ if 
$(2\pi \hk)\cdot\o_g>0$, or dually, if $k\cdot H>0$, and $\ssw_{X,H}(k) =\ssw_{X,H}^-(k)$ if $k\cdot H <0$. As in the case above, we work with the modified invariant 
$\swu_{X,H}(k)=\sum\limits_{c(\fs)=k}\ssw_{X,H}(\fs)$.

We now proceed to the proof of the theorem.

\begin{proof}[Proof of Theorem~\ref{SW}] Set $X_0=X_{\L,\lam}(0)$. Recall that $\L_0$ is the torus in $X_0$ which is the core torus of the surgery. There is a surface in $X_0$ which intersects $\L_0$ once, and it follows that $\L_0$ is essential (in fact primitive in $H_2$). The surface in question is the union of a normal disk to $\L_0$ and the surface in $X_0\- N_{\L_0} = X\- N_{\L}$ bounded by $S^1_{\lam}$. 

Let $X_n=X_{\L,\lam}(1/n)$. Then $\L_{1/n}$ is the core torus of the surgery in $X_n$. Its meridian $\mu_{\L_{1/n}}$ represents $n[\lam]+[\mu_{\L}]$, which in 
$X_n\- N_{\L_{1/n}} = X\- N_{\L}$ is homologous to $[\mu_{\L}]$, a nontorsion class. This means that $\L_{1/n}$ is nullhomologous in $X_n$.

Let $k_n\in  H_2(X_n;\Z)$ be characteristic. The classes in $H_2(X_0;\Z)$ which are characteristic and which agree with the restriction of $k_n$ in $H_2(X_n\- N_{\L_{1/n}},\bd;\Z)=H_2(X\- N_{\L},\bd;\Z)$ in the diagram:
\[ \begin{array}{ccc}
H_2(X_n;\Z) &\longrightarrow & H_2(X_n, N_{\L_{1/n}};\Z)\\
&&\Big\downarrow \cong\\
&&H_2(X\- N_{\L},\bd;\Z)\\
&&\Big\uparrow \cong\\
H_2(X_0;\Z)&\longrightarrow & H_2(X_0,N_{\L_0};\Z)
\end{array}\]
differ by an even multiple of $[\L_0]$. The adjunction inequality implies that each such class with a nontrivial Seiberg-Witten invariant is orthogonal to $[\L_0]$. There is just one class in $X_n$ which restricts to a fixed class in  $H_2(X_n\- N_{\L_{1/n}},\bd;\Z)$ because $\L_{1/n}$ is nullhomologous, and the same is true for $X$ (and $\L$).

It follows from \cite{MMS} that 
\[ \swu_{X_n}(k_n) = \swu_X(k) + n\sum_i \swu_{X_0}(k_0+2i[\L_0]) \]
where $k$ and $k_0$ are classes which restrict as in the paragraph above. Recalling the hypothesis that there is a $k_0\in H_2(X_0;\Z)$, satisfying $\sum\limits_i\swu_{X_0}(k_0+2i[\L_0])\ne0$,
it follows that the integer invariants
\[ S_n = \max\{ |\swu_{X_n}(k_n)|;\; k_n \, {\text{basic for $X_n$}} \} \]
will distinguish an infinite family of pairwise nondiffeomorphic manifolds among the $X_n$.

In case $b^+_X=1$, we need to check issues with chambers. The inclusions of $X\- N_{\L}$ in $X$ and $X_n$ induce isomorphisms on $H_2$ and thus an isomorphism of $H_2(X_n;\Z)$ with $H_2(X;\Z)$. The gluing formula of \cite{MMS} relates chambers using this isomorphism. So, for example, if $\ssw_{X,H}(k)= \ssw^+_{X,H}(k)$ this means that $k\cdot H>0$. The isomorphism above gives a $k_n\in H_2(X_n;\Z)$ and an $H\in H_2(X_n;\R)$ (and $H_n\cdot H_n=H\cdot H>0$), and it also gives $k_n\cdot H_n=k\cdot H>0$; so $\ssw_{X_n,H_n}(k_n)=\ssw^+_{X_n,H_n}(k_n)$. Thus the gluing formula applies to the invariant $\ssw_{X,H}$. This works for any choice of period point $H$. Hence the argument in the $b^+>1$ case applies directly to $b^+=1$ as well.
\end{proof}

Corollary \ref{CSW} follows similarly. When $X_0$ has just one basic class up to sign, it is straightforward to see that  if $n\ne m$ the collection of values of the Seiberg-Witten invariant $\{\swu_{X_n}(k)\}$ is different from $\{\swu_{X_m}(k)\}$.

\begin{cor}\label{SWcor} Let $X$ be a smooth closed oriented $4$-manifold which contains a nullhomologous torus $\L$ and let $\lam$ be a simple loop on $\L$ so that $S^1_{\lam}$ is nullhomologous in $X\- N_{\lam}$. Suppose also that there is a square $0$ torus $T\C X_0$ that satisfies $T\cdot\L_0\ne0$. If $X_0$ has a basic class, {\it{i.e.}} a class $k_0$ with $\ssw'_{X_0}(k_0)\ne0$,
then among the manifolds $\{ X_{\L,\lam}(1/n)\}$, infinitely many are pairwise nondiffeomorphic.
\end{cor}
\begin{proof} If $T\C X_0$ is a torus of square $0$ that satisfies $T\cdot\L_0\ne0$, then the adjunction inequality implies that in each collection $\{k_0+2i\L_0\}$, there is at most one basic class. Thus the hypothesis of Theorem~\ref{SW} is satisfied.
\end{proof}

We now provide two examples to illustrate the reverse engineering procedure.

\section{Fake $\CPC$'s}

The 2-fold symmetric product $Y=Sym^2(\Sig_3)$ is the quotient of $\Sig_3 \x \Sig_3$ by the action of the involution $t:\Sig_3 \x \Sig_3 \to \Sig_3 \x \Sig_3$ given by $t(x,y)=(y,x)$. Let $\{a_i, b_i\}$, $i=1,2,3$
denote standard generators for $\pi_1(\Sig_3)$.  
It follows from \cite{Per} that  the natural singular  K\"ahler form $Sym^2(\omega)$ on $Y$ derived from the K\"ahler curve $(\Sig_3,\omega)$ admits a cohomologous smoothing to a K\"ahler form which equals $Sym^2(\omega)$ away from a chosen neighborhood of the diagonal. (We thank Paul Kirk for pointing out the necessity of this reference.) 

We obtain a basis for $H_{2}(Y;\Z)\cong \Z^{16}$ as follows.  The tori $a_{i}\x a_{j}$, $b_{i}\x b_{j}$, $a_{i} \x b_{j}$, and $b_i\x a_j$, $i < j$, in $\Sig_3 \x \Sig_3$ descend to twelve tori  of self-intersection $0$ in $Y$, and we also denote these by $a_{i}\x a_{j}$, $b_{i}\x b_{j}$, $a_{i} \x b_{j}$, and $b_{i} \x a_{j}$. The three tori $a_{i} \x b_{i}$ in $\Sig_3 \x \Sig_3$ descend to  tori $T_i$ of square $-1$, and together with the image of $\{\text{pt}\} \x \Sig_3 \cup \Sig_3 \x \{\text{pt}\}$, a genus 3 surface which represents a homology class $b$ with self-intersection $+1$, we get a basis for $H_{2}(Y;\Z)$. The euler number $\eu(Y)=6$ and its signature $\sign(Y)=-2$, in agreement with the characteristic numbers for $\CPC$.

\begin{figure}[!ht]
\begin{center}
\includegraphics[scale=1.15]{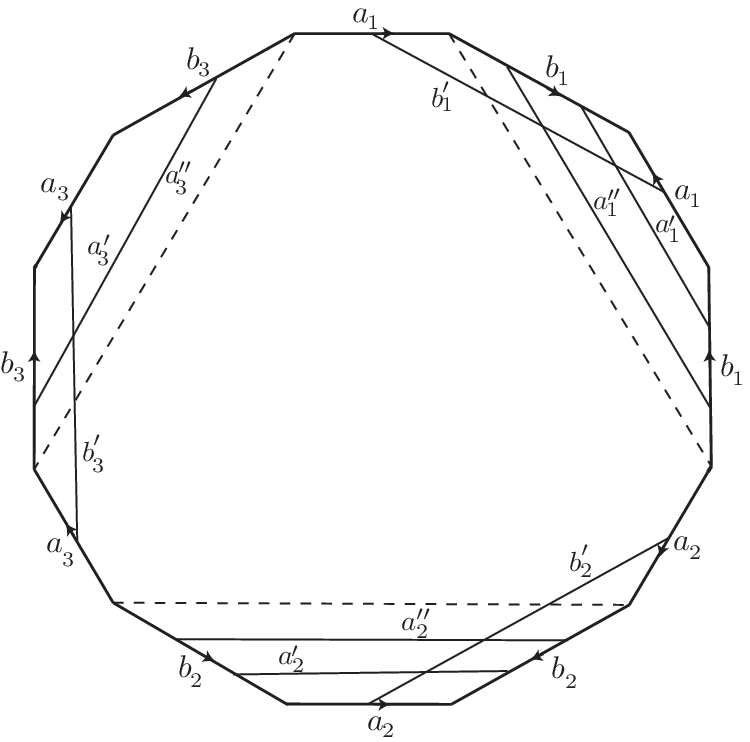}\vspace{-.2in}
\end{center}
\caption{}
\label{Figure1}
\end{figure}

To establish some notation, consider Figure 1. For example, we see loops $a_i$, $a_i'$, and $a_i''$.  We also have based loops (with basepoint $x$, the vertex) which we shall denote by $\a_i$, $\a_i'$, $\b_j$, etc. The based loop $\a_2'$, for example, is the one which starts at the lower left vertex $x$, proceeds backwards along $b_2$ to the initial point $a_2'(0)$ of $a_2'$, then traverses $a_2'$ until it gets to its endpoint $a_2'(1)$, and then heads vertically downward back to $x$. The based loop $\a_2''$ starts at the vertex at the initial point of $b_2$, travels along $b_2$ to $a_2''(0)$, traverses $a_2''$, then at $a_2''(1)$ it heads upward back to $x$. The based loop $\a_i$ is equal to $a_i$.

Using $(x,x)$ as basepoint, the abelian group $\pi_{1}(Y)=\Z^{6}$ is generated by the $\a_i=\a_{i}\x \{x\}$ and $\b_j=\b_{j}\x \{x\}$. We will perform six surgeries on disjoint Lagrangian tori to kill these generators. We need some notation to describe our surgeries. First notice that the Lagrangian tori $a_1'\x a_2'$, $a_1''\x b_2'$,
$a_1'\x a_3'$, $b_1'\x a_3''$, $a_2'\x a_3'$, and $a_2''\x b_3'$ in $\Sig_3 \x \Sig_3$ are disjoint from the diagonal, and so they descend to Lagrangian tori in $Y=Sym^2(\Sig_3)$. (We still denote these tori in $Y$ by $a_1'\x a_2'$, etc.) They are mutually disjoint.
If, for example, we consider the torus $a_1'\x a_2'$ and do $n$-framed surgery in $Y$ along the loop $a_1'$ with respect to the Lagrangian framing, we denote this as the surgery $(a_1'\x a_2', a_1',n)$.  We now perform surgeries along disjoint Lagrangian tori 
\begin{gather*}  (a_1'\x a_2', a_2',-1),\ \ (a_1''\x b_2', b_2',-1), \ \ (a_1'\x a_3', a_1',-1), \\
(b_1'\x a_3'', b_1',-1), \ \ (a_2'\x a_3',a_3',-1), \ \ (a_2''\x b_3', b_3',-1).
\end{gather*}

Denote these six Lagrangian tori by $\L_k$, $k=1,\dots,6$. We consider tubular neighborhoods $N_{\L_k}$ with their Lagrangian framings. For example, if $\L_1= a_1'\x a_2'$ then $N_{\L_1}$ can be visualized in Figure~1 as the product of the two annuli given by the trapezoids between $a_i'$ and $\ba_i'$, $i=1,2$. Then we get Lagrangian pushoffs $\ha_1'=
\ba_1'\x\{\ba_2'(0)\}$ and $\ha_2'=\{\ba'_1(0)\}\x\ba_2'$ on $\bd N_{\L_1}$. Let 
$m_{\L_1}$ be the boundary circle to the normal disk passing through the intersection point 
$x_{\L_1}=\{\ba_1'(0)\}\x \{\ba_2'(0)\}$ of $\ha_1'$ and $\ha_2'$.
In like manner, for each $k$ we have a `distinguished triple of loops' $\{\hg_{k,1},\hg_{k,2};m_{\L_k}\}$ which pass through $x_{\L_k}$ on  $\bd N_{\L_k}$.

\begin{lem} For each $\L_k$ there is a path, $\g_k$, from the basepoint $(x,x)$ of $Y$ to 
$x_{\L_k}$ which except for its endpoint lies in $X\- \bigcup\limits_1^6 N_{\L_j}$. 
If we base the distingushed triple of loops for $\L_k$ at $x$ using $\g_k$ and $\g_k^{-1}$, then the corresponding based loops are respectively homotopic in $X\- \bigcup\limits_1^6 \L_j$ to 
\begin{eqnarray*}
&& \{\a_1, \a_2; [\b_1^{-1},\b_2^{-1}]\}, \quad 
\{ \b_1\a_1\b_1^{-1}, \b_2; [\b_1,\a_2^{-1}]\},\\
&& \{ \a_1,\a_3; [\b_1^{-1},\b_3^{-1}]\}, \quad
\{\b_1, \b_3\a_3\b_3^{-1}; [\a_1^{-1},\b_3]\},\\
&& \{\a_2, \a_3; [\b_2^{-1},\b_3^{-1}]\},\quad 
\{\b_2\a_2\b_2^{-1}, \b_3; [\b_2,\a_3^{-1}]\}.
\end{eqnarray*}
\end{lem}

\begin{proof}(cf. \cite{BK}) We give the proof for the first two tori $\L_1=a_1'\x a_2'$ and 
$\L_2=a_1''\x b_2'$. The other cases are similar. First consider $\L_1$. Let $\z_1(t)$ be the path from the basepoint $x$ to $\ba_1'(0)$ traveling backwards along $b_1$ as discussed above, and let $\eta_1(t)$ be a similar path starting at $x$ and ending at $\ba_2'(0)$, traveling backwards along $b_2$. Set $\g_1(t) = \z_1(t)\x\eta_1(t)$.
The trapezoid in Figure~1 between $\ba_1'$ and $a_1$ gives rise to a homotopy between the based loop $\g_1\cdot\ha_1'\cdot\g_1^{-1}$ and $\a_1=a_1$.

For each $s\in [0,1]$, let $\z_{1,s}(t)=\z_1(st)$ and $\eta_{1,s}(t)=\eta_1(st)$.
Thought of as a homotopy, the trapezoid between $\ba_1'$ and $a_1$ is composed of parallel closed paths $\ba'_{1,s}$, where $\ba'_{1,0}=\ba_1'$, $\ba'_{1,1}=a_1$, and the initial point of $\ba'_{1,s}$ is $\z_1(1-s)$. This defines the based homotopy of $\g_1\cdot\ha_1'\cdot\g_1^{-1}$ to $\a_1$ whose path at level $s$ is the image in $Y$ of the product of paths
\[ (\z_{1,1-s}\x\eta_{1,1-s})\cdot(\ba'_{1,s}\x \{\eta_1(1-s)\})\cdot (\z_{1,1-s}^{-1}\x \eta_{1,1-s}^{-1}).\]
It is easily seen that the trace of this homotopy is disjoint from all the tori $\L_k$. Similarly, the trapezoid between $\ba_2'$ and $a_2$ gives the homotopy between the based loop $\g_1\cdot\ha_2'\cdot\g_1^{-1}$ and $\a_2=a_2$. 

The meridian $m_{\L_1}$ lies in the orthogonal torus $b_1\x b_2$ which is Lagrangian away from the basepoint $(x,x)$ of $Y$. This torus also contains the base path $\g_1$. Since $\g_1$ runs backwards along both $b_1$ and $b_2$ to $x_{\L_1}$, it follows that when $m_{\L_1}$ is based using $\g_1$, it is given by the commutator of the appropriately oriented $\pi_1$-generators of $b_1\x b_2$; {\it{viz.}} $[\b_1^{-1},\b_2^{-1}]$. (There is a choice of orientation for this meridian. The opposite choice would cause us to change the signs of our surgeries. For definiteness, we choose the orientation implied by the statement of the lemma.)

Next consider $\L_2 = a_1''\x b_2'$. Its tubular neighborhood $N_{\L_2}$ with Lagrangian framing is given by the product of the two annuli described by trapezoids in Figure~1 between $a_1''$ and $\ba_1''$ and between $b_2'$ and $\bar{b}_2'$. (The curve $\bar{b}_2'$ is not shown in Figure~1, but it is analogous to $\ba_1'$.)
Let $\z_2(t)$ be the path traveling positively along $b_1$ from $x=b_1(0)$ to $\ba_1''(0)$, and let $\eta_2(t)$ be the path traveling backwards along $a_2$ from $x$ to $\bar{b}_2'(0)$. Set $\g_2(t) = \z_2(t)\x\eta_2(t)$. Then, analogously to the argument above, the trapezoid between $\ba_1''$ and the nearby parallel dashed line describes a homotopy from $\g_2\cdot \ha_1''\cdot\g_2^{-1}$ to the based loop given by this dashed line --- which is $\b_1\a_1\b_1^{-1}$. The homotopy from $\g_2\cdot\hat{b}_2'\cdot\g_2^{-1}$ to $\b_2$ is given by the trapezoid between $\bar{b}_2'$ and $b_2$. To see the meridian $m_{\L_2}$ we use the dual torus $b_1\x a_2$ which contains the path $\g_2$
which runs positively along $b_1$ and negatively along $a_2$, and the same reasoning as above shows that when $m_{\L_2}$ is based using $\g_2$, it is given by $ [\b_1,\a_2^{-1}]$. 
\end{proof}

Denote the manifold obtained from these six surgeries by $X$. The result of each surgery is to reduce the first betti number by one, reduce the second betti number by two, and introduce a relation in $\pi_1$. For example, because of Lemma~1, the surgery  $(a_1'\x a_2', a_2',-1)$ introduces the relation $\a_2= [\b_1^{-1},\b_2^{-1}]$, and the surgery $(a_1''\x b_2', b_2',-1)$ introduces the relation $\b_2= [\b_1,\a_2^{-1}]$. 
Note that $b_1(X)=0$, and since the surgeries change neither the euler number nor signature, $b_2(X)=4$, and $X$ is a rational homology $\CPC$.

The following relations hold in $\pi_1(X)$:
\begin{gather*} 
\a_2= [\b_1^{-1},\b_2^{-1}],\ \  \b_2= [\b_1,\a_2^{-1}],\ \ \a_1= [\b_1^{-1},\b_3^{-1}],\\
\b_1= [\a_1^{-1},\b_3],\ \ \a_3= [\b_2^{-1},\b_3^{-1}],\ \  \b_3=[\b_2,\a_3^{-1}],\\
[\a_1,\b_1]=1, \ \ [\a_1,\a_2]=1,\ \  [\a_1,\b_2]=1,\ \ [\a_1,\a_3]=1,\ \ [\b_1,\a_3]=1,\\ 
 [\a_2,\b_2]=1,\ \  [\a_2,\a_3]=1,\ \ [\a_2,\b_3]=1,\ \ [\a_3,\b_3]=1.
\end {gather*} 
Thus we have $\b_2= [\b_1,\a_2^{-1}]=[[\a_1^{-1},\b_3],\a_2^{-1}]=1$, using the commutativity relations $[\a_2,\b_3]=1$ and $[\a_1,\a_2]=1$. Now it follows from the other relations that $\pi_1(X)=1$. 

Since the surgeries that we perform on the Lagrangian tori all have surgery coefficient $\pm1$ with respect to the Lagrangian framing, the resultant manifolds all have induced symplectic structures. One simple way to see that $X$ is not diffeomorphic to $\CPC$ is to use the fact from \cite{LL} that $\CPC$ has a unique symplectic form up to diffeomorphism and symplectic deformation. This means that for any symplectic form on $\CPC$, the canonical class must pair negatively with the symplectic form. On $Y=Sym^2(\Sig_3)$, which is a surface of general type, the canonical class pairs positively with the symplectic form, and since we have constructed $X$ by surgeries on Lagrangian tori of $Y$, the same is still true in $X$. (The point here is that if $\hat{Y}$ is the result of a Luttinger surgery on $Y$, then the complements of tubular neighborhoods of the respective Lagrangian tori in each can be identified, and the restrictions of the symplectic forms can as well. The canonical classes are supported in the complements of these tori and agree over the complements of the tubular neighborhoods. It follows that as elements of $H_2(\hat{Y};\Z)$ the Poincar\'e duals satisfy $\hat{K}=K+n[\hat{T}]$. Cf. \cite{ADK}. So $\hat{K}\cdot \hat{\o} = K\cdot\o$ since $\hat{T}$ is Lagrangian.) Hence $X$ cannot be diffeomorphic to $\CPC$. 

\begin{thm}[cf. \cite{AP,BK}] The symplectic manifold $X$ is irreducible and homeomorphic but not diffeomorphic to $\CPC$.
\end{thm}

The irreducibility of $X$ follows from \cite{HK} once we show that $X$ is minimal. This follows from the Seiberg-Witten calculations below.
It is interesting to ask whether $X$ is actually diffeomorphic to the symplectic manifolds constructed in \cite{AP,BK}.

In order to produce an infinite family of exotic $\CPC$'s, let $X'$ denote the result of the first five Luttinger surgeries on $Y$. Thus $b_1(X')=1$ and $b_2(X')=6$. We construct $X$ by performing a surgery $(a_2''\x b_3', b_3',-1)$ in $X'$. In $X$, the surgery gives us a nullhomologous torus  $\L$, the ``core" of the surgery. There is a loop $\lam$ on $\L$ so that surgery on $(\L,\lam)$ gives $X'$ back. The framing for this surgery must be the nullhomologous framing.  We apply Theorem~\ref{SW} to $(X,\L,\lam)$.  In fact, Corollary~\ref{CSW} will tell us that the manifolds $X_{\L,\lam}(1/n)$
are pairwise nondiffeomorphic once we see that $X'$ has exactly two basic classes.  (The manifold denoted $X_0$ in Corollary~\ref{CSW} is our manifold $X'$.)  Note that $X_n = X_{\L,\lam}(1/n)$ is the result of performing the surgery $(a_2''\x b_3', b_3', n+1)$ in $X'$. (These are {\em{not}} Luttinger surgeries.)

\begin{thm}\label{inf} The manifolds $\CPC$, $X$, and $X_n$, $n\ge2$, are pairwise homeomorphic and (except for $\CPC$) are minimal, but no two are diffeomorphic.
\end{thm}

\begin{proof} The homeomorphism statement will follow once we see that each $X_n$ is simply connected. A presentation for $\pi_1(X_n)$ is obtained from the one above for $ \pi_1(X)$ by replacing the relation $\b_3=[\b_2,\a_3^{-1}]$ by 
$\b_3=[\b_2,\a_3^{-1}]^{-(n+1)}$, and $\pi_1(X_n)=1$ follows as above.

Next, we need to show that the manifold $X'$ has just two basic classes, $\pm$ its canonical class, and then call on Corollary~\ref{CSW}. Since $Y$ is a surface of general type, its only basic classes are $\pm$ its canonical class, i.e. $3b+T_1+T_2+T_3$, where $b$ and $T_i$ are the classes described above.
According to \cite{MMS},  each time we do a surgery, the Seiberg-Witten invariant of the result is calculated in terms of the Seiberg-Witten invariants of the original manifold and those of the result of the surgery that kills the curve on the torus. For example, if $Y_1$ is the result of the surgery $(a_1'\x a_2', a_2',-1)$ on $Y$, then let $Z$ be the result of the surgery that kills $a_2'$ directly ($0$-surgery). In $Z$, the surface $\Sig_3 \x \{\text{pt}\}$, which represents $b$, has its genus reduced by one. Applying the adjunction inequality to this situation, we see that any basic class of $Z$ has the form $\pm  b\pm T_1\pm T_2 \pm T_3$. Since the square of a basic class must be $3\,\sign(Z) + 2\,\eu(Z) = 6$, in fact none of these classes can be basic; so the Seiberg-Witten invariant of $Z$ vanishes.
The result of this argument is that the manifold $Y_1$ also has just two basic classes, $\pm$ its canonical class.
The very same argument works for each surgery and finally shows that $X'$ has just two basic classes.

Thus $X$ and $X_n$ have just two basic classes, $\pm k_n$, and the difference is a class of square $(2k_n)^2= 24$. If one of these manifolds failed to be minimal, it would have to have a pair of basic classes, $k\pm E$, whose difference has square $-4$. Thus $X$ and $X_n$ are minimal.
\end{proof}

In order to obtain infinitely many smooth structures, we did not need to perform this last step which shows that $X$ and $X_n$ have just two basic classes. We did this to explicitly show that all the $X_n$ are distinct. The hypothesis of Corollary~\ref{SWcor} is satisfied because each of the Lagrangian tori on which surgery is performed has a dual torus of square $0$.

Each $X_{n}$ contains disjoint embeddings of a minimal genus $3$ surface representing $b$ and the three tori $T_{1}$, $T_{2}$, $T_{3}$ with self-intersection $-1$. More interestingly, each $X_{n}$ contains a sphere of self-intersection $-2$ representing $b -[T_{1}] -[T_{2}] -[T_{3}]$ that is the image in $Sym^2(\Sig_3)$ of a pushoff of the diagonal in $\Sig_3 \x \Sig_3$. These surfaces can be useful for other constructions.

The symmetric product $Z_{\ell}=Sym^2(\Sig_{\ell})$ of a genus $\ell$\/ surface with itself  has $\pi_1(Sym^2(\Sig_{\ell}))=H_1(\Sig_{\ell};\Z)$ and  that 
$\eu(Z_{\ell})=(\ell-1)(2\ell-3)=2\ell^2-5\ell+3$, $\sign(Z_{\ell})=1-\ell$, and $b_2(Z_{\ell})=2\ell^2-\ell+1$. Thus $Z_{\ell}$ is a model for $(\ell^2-3\ell+1)\CP\#(\ell^2-2\ell)\CPb$. A straightforward generalization of the above application of reverse engineering provides infinitely many distinct smooth structures on these manifolds, one of which is  symplectic.

\section{Fake homology $\SS$'s}

We now give an example to point out that the computation of fundamental groups in the reverse engineering procedure can be difficult.

\begin{figure}[ht]
\begin{center}
\includegraphics[scale=.45]{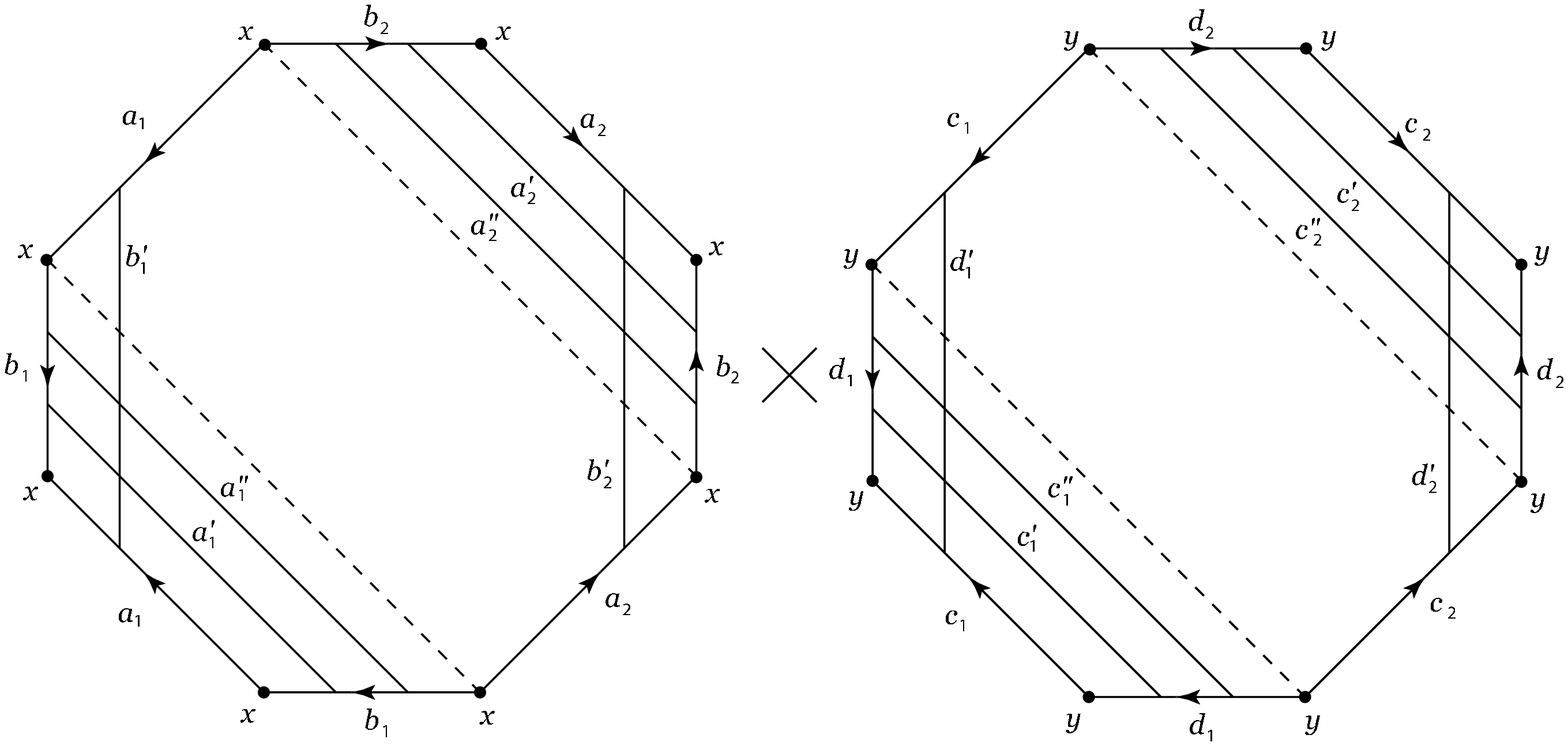}
\end{center}
\caption{}
\label{fig:Lagrangian-pushoff}
\end{figure}

Let $Y=\Sig_2\x\Sig_2$, the product of two genus $2$ surfaces, and denote the standard generators of $\pi_1$ by $\{ a_i,b_i\}$ and $\{c_i,d_i\}$ for $i=1,2$. (Since the calculations in this section parallel those of \S3, for the sake of simplicity we no longer denote the distinction between based and unbased loops.) So
$\pi_1(Y)$ has these eight generators with relations $[a_1,b_1][a_2,b_2]=1$, $[c_1,d_1][c_2,d_2]=1$ and all $a_i$ and $b_i$ commute with all $c_j$ and $d_j$. The area forms on the two copies of $\Sig_2$ induce a symplectic form on the product $Y$, all tori of the form $a_i\x c_j$,  $a_i\x d_j$, $b_i\x c_j$, and $b_i\x d_j$ are Lagrangian, and the Lagrangian framing is the obvious one coming from the product structure. The euler number $\eu(Y)=4$ and its signature $\sign(Y)=0$, in agreement with the characteristic numbers for $\SS$.

Perform eight Luttinger surgeries along the Lagrangian tori 
\begin{eqnarray*}
&&(a_1' \times c_1', a_1', -1), \quad (b_1' \times c_1'', b_1', -1), \quad
(a_2' \times c_2', a_2', -1), \quad (b_2' \times c_2'', b_2', -1),\\
&&(a_2' \times c_1', c_1', +1), \quad (a_2'' \times d_1', d_1', +1),\quad
(a_1' \times c_2', c_2', +1), \quad (a_1'' \times d_2', d_2', +1).
\end{eqnarray*} 
to obtain a symplectic manifold $X$. (See Figure~2.)

As in \S3 each of these Lagrangian tori have tubular neighborhoods on which there is a distinguished triple of curves and arguing as in $\S3$ we have: 

\begin{lem} Let $Y'$ be the complement of the above eight Lagrangian tori in 
$Y=\Sigma_2\times\Sigma_2$.  Inside $Y'$, there are basepaths from the basepoint of $Y$ to the basepoint of the boundaries of the tubular neighborhoods of these tori, so that if we base the distinguished triple of loops using these basepaths 
 then the corresponding based loops are homotopic in $Y'$ to 
\begin{eqnarray*}
&& \{ a_1,c_1; [b_1^{-1},d_1^{-1}]\}, \quad 
\{ b_1, d_1c_1d_1^{-1}; [a_1^{-1},d_1]\},\\
&& \{ a_2, c_2; [b_2^{-1},d_2^{-1}]\}, \quad
\{ b_2, d_2c_2d_2^{-1};[a_2^{-1},d_2]\} ,\\
&& \{ a_2, c_1; [b_2^{-1},d_1^{-1}]\},\quad 
\{ b_2a_2b_2^{-1}, d_1; [b_2,c_1^{-1}]\},\\
&& \{ a_1, c_2; [b_1^{-1},d_2^{-1}]\},\quad 
\{ b_1a_1b_1^{-1}, d_2; [b_1,c_2^{-1}]\}.
\end{eqnarray*}
\end{lem}

The following relations hold in $\pi_1(X)$:
\begin{gather*} 
[b_1^{-1},d_1^{-1}]=a_1,\ \  [a_1^{-1},d_1]=b_1,\ \  [b_2^{-1},d_2^{-1}]=a_2,\ \  [a_2^{-1},d_2]=b_2,\\
[d_1^{-1},b_2^{-1}]=c_1,\ \ [c_1^{-1},b_2]=d_1,\ \ [d^{-1}_2,b^{-1}_1]=c_2,\ \ [c_2^{-1},b_1]=d_2,\\
 [a_1,c_1]=1, \ \ [a_1,c_2]=1,\ \  [a_1,d_2]=1,\ \ [b_1,c_1]=1,\\ 
[a_2,c_1]=1, \ \ [a_2,c_2]=1,\ \  [a_2,d_1]=1,\ \ [b_2,c_2]=1,\\
[a_1,b_1][a_2,b_2]=1,\ \ [c_1,d_1][c_2,d_2]=1.
\end {gather*} 

From these relations it is clear that $H_1(X:\Z)=0$, and since the surgeries change neither the euler number nor signature, $b_2=2$. In fact, the only homology classes that survive are those represented by $\Sig_2\x\{\text{pt}\}$ and $\{\text{pt}\}\x \Sig_2$; so $X$ is a homology $\SS$.

We have been unable to determine if the perfect group $\pi_1(X)$ is trivial or not. Also, there are other surgeries that can be performed, and also other sets of eight Lagrangian tori that can be surgered, to obtain many presentations of perfect groups that we have not succeeded in showing are trivial.

We can produce an infinite family of distinct homology $\SS$'s, exactly as in \S3. 
The presentation for $\pi_1$ is exactly the one given above for $\pi_1(X)$ except that the relation $[c_2^{-1},b_1]=d_2$ is replaced by 
$[c_2^{-1},b_1]^{n+1}=d_2$ compounding the difficulty of determining whether the group is trivial.

\end{document}